\documentclass[a4paper,11pt]{amsart}
\usepackage[utf8]{inputenc}
\usepackage[T1]{fontenc}
\usepackage[numbers]{natbib}
\usepackage[english]{babel}
\usepackage[margin=3cm]{geometry}
\usepackage{amssymb}
\usepackage{dsfont}
\usepackage{enumerate}
\usepackage[inline]{enumitem}
\usepackage{mathtools}
\usepackage{amsmath}
\usepackage{color}
\usepackage{bbm}         
\usepackage{amsfonts}
\usepackage{graphicx,verbatim}
\usepackage{amsthm}
\usepackage{url}
\usepackage{tikz,multicol,manfnt}
\usepackage{qtree}
\usepackage[weather,clock,alpine]{ifsym}
\usepackage{xparse}
\usepackage{parnotes}
\usepackage{placeins}
\usepackage{xspace}
\usepackage{xcolor}
\usepackage{tikz-cd}
\usetikzlibrary{positioning}
\usetikzlibrary{decorations.pathreplacing}
\usetikzlibrary{calc,arrows,shapes,fit}

\usepackage{hyperref}
\usepackage{cleveref}

\newcommand{\rond}{\mathcal}

\newcommand{\satisf}{\vDash}
\newcommand{\et}{\wedge}

\newcommand{\Ou}{\bigvee}
\renewcommand{\phi}{\varphi}
\renewcommand{\epsilon}{\varepsilon}

\NewDocumentCommand{\set}{mg}{\left\{#1\IfNoValueF{#2}{\;\middle\vert\;#2}\right\}}

\newcommand{\IP}[1]{IP$_{\!#1}$\xspace}
\newcommand{\NIP}[1]{N\IP{#1}}
\newcommand{\IPn}{\IP n}
\newcommand{\NIPn}{\NIP n}
\newcommand{\tpp}{\!\cdots\!}
\newcommand{\ignore}[1]{}

\newcommand{\subtitle}[1]{%
  \posttitle{%
    \par\end{center}
    \begin{center}\large#1\end{center}
    \vskip0.5em}%
}
\newcommand{\subtitles}[3]{%
  \posttitle{%
    \par\end{center}
    \begin{center}\large#1\\\normalsize#2\\\small#3\end{center}
    \vskip0.5em}%
}

\DeclareMathOperator{\ch}{char}

\DeclareMathOperator{\qftp}{qftp}
\DeclareMathOperator{\tp}{tp}

\theoremstyle{plain}

\newtheorem{thm}{Theorem}[section]
\newtheorem*{thm*}{Theorem}
\newtheorem{cor}[thm]{Corollary}
\newtheorem*{cor*}{Corollary}
\newtheorem{lem}[thm]{Lemma}
\newtheorem{prop}[thm]{Proposition}
\newtheorem*{prop*}{Proposition}
\newtheorem{con}[thm]{Conjecture}
\newtheorem*{claim}{Claim}

\theoremstyle{definition}
\newtheorem{dfn}[thm]{Definition}
\theoremstyle{remark}
\newtheorem{rem}[thm]{Remark}
\newtheorem{ex}[thm]{Example}

\usepackage{etoolbox}
\makeatletter
\patchcmd{\@maketitle}{\normalsize}{\LARGE}{}{}
\makeatother
 
  \author{Blaise Boissonneau}
  \address{Blaise Boissonneau\\Heinrich-Heine-Universität Düsseldorf, 40225 Düsseldorf, Germany}
\email{blaise.boissonneau@hhu.de}
\thanks{The author was funded by Franziska Jahnke's fellowship from the Daimler and Benz foundation. This research was also partially funded by the DAAD through the `Kurzstipendien für Doktoranden 2020/21', the MSRI via the Decidability, Definablility and Computability programme, the KU Leuven, and the HHU Düsseldorf.}

  \title{\texorpdfstring{\NIPn}{NIPn} CHIPS}
  \date{\today}

\begin{document}

  \begin{abstract}
  We give general conditions under which classes of valued fields have \NIPn transfer and generalize Anscombe-Jahnke's classification of NIP henselian valued fields to \NIPn henselian valued fields.
  \end{abstract}
  
  \maketitle

  \section{Introduction}
  
  The main goal of this paper is to prove the following classification of \NIPn henselian valued fields:
  
  \begin{thm}\label{NIPthm}
 Let $(K,v)$ be a henselian valued field. Then $(K,v)$ is \NIPn iff the following holds:
\begin{enumerate}
\item the residue field $k$ is \NIPn, and
\item either
\begin{enumerate}[ref=\theenumi{}.(\alph*)]
\item\label{forbelow} $(K,v)$ is of equicharacteristic and is either trivial or SAMK, or
\item $(K,v)$ has mixed characteristic $(0,p)$, $(K,v_p)$ is finitely ramified, and $(k_p,\overline{v})$ satisfies \ref{forbelow}, or
\item $(K,v)$ has mixed characteristic $(0,p)$ and $(k_0,\overline{v})$ is AMK.
\end{enumerate}
\end{enumerate}
\end{thm} 

 In this paper we write ``(S)AMK'' for ``(separably) algebraically maximal Kaplansky'', and we denote by $v_0$ the finest valuation of residue characteristic 0 and $v_p$ the coarsest valuation of residue characteristic $p$; we refer to \cref{subsec-not} for details.
  
  For $n=1$, that is, for NIP henselian valued fields, \Cref{NIPthm} is a theorem of Anscombe and Jahnke, see \cite[Thm.~5.1]{AJ-NIP}. For the most part, in order to extend their classification to arbitrary $n$, we follow their strategy, mutatis mutandis.
  
  In a precedent article, we proved the following:
  
  \begin{thm}[{\cite[Thm.~1.2]{itme-asext}}]\label{asext}
   Let $(K,v)$ be a $p$-henselian valued field. If $K$ is \NIPn, then either:
\begin{enumerate}[label=(\alph*)]
\item\label{eqchar} $(K,v)$ is of equicharacteristic and is either trivial or SAMK, or
\item $(K,v)$ has mixed characteristic $(0,p)$, $(K,v_p)$ is finitely ramified, and $(k_p,\overline{v})$ satisfies \ref{eqchar}, or
\item $(K,v)$ has mixed characteristic $(0,p)$ and $(k_0,\overline{v})$ is AMK.
\end{enumerate}
  \end{thm}

  Thus, one direction of \Cref{NIPthm} is already proven; it is even slightly stronger as we work with $p$-henselian valuation, and only the pure field $K$ is assumed to be \NIPn.
  
  It remains to prove the other direction, which is a \NIPn transfer theorem. A transfer theorem is of the form ``if the residue field and the value group of $(K,v)$ are (*), then $(K,v)$ is (*)'', (*) can be any suitable model-theoretic condition; in this paper, we will talk about NIP, \NIPn, and NTP2 transfer theorems, and these theorems usually assume strong conditions on $(K,v)$, such as henselianity or algebraic maximality.
  
  \subsection{A short history of transfer theorems}
  
  Transfer theorems have been studied as early as 1981, with the following result of Delon:
  
  \begin{thm}[\cite{Delon}]
  Let $(K,v)$ be henselian of equicharacteristic 0, then $(K,v)$ is NIP iff its residue field $k$ and its value group $\Gamma$ are NIP.
  \end{thm}
  
  Note that the condition on the value group is empty since Gurevich and Schmidtt proved in \cite{GS} that all pure ordered abelian groups are NIP, therefore, more recent NIP transfer theorems don't include this clause.
  
  In other characteristics, more assumptions are needed:
  
  \begin{thm}[\cite{Bel}]\ 
  \begin{enumerate}
   \item Let $(K,v)$ be henselian, of equicharacteristic $p$, and AMK; then $(K,v)$ is NIP iff its residue field $k$ is NIP.
   \item Let $(K,v)$ be henselian, of mixed characteristic, unramified, and with residue field perfect; then $(K,v)$ is NIP iff its residue field $k$ is NIP.
   \end{enumerate}
  \end{thm}
  
  All of these results were recently generalized by Anscombe, Jahnke and Simon by using a powerful method, developed first in 2014 by Chernikov and Hils for NTP2 transfer. Given a complete theory of valued fields, possibly with augmented structure, consider the following properties:
  
  \begin{enumerate}
  \item[(Im):] For any small model $K$ and any singleton $b$ (from a monster model) such that $K(b)/K$ is immediate, we have that $\tp(b/K)$ is implied by instances of NTP2 formulas, that is, there is $p\subset\tp(b/K)$ closed under conjunctions and such that:
   \begin{itemize}
    \item any formula $\phi(x,y)\in p$ is NTP2,
    \item for any formula $\psi(x,y)$, $\psi(b,K)$ holds iff $p\vdash\psi(x,K)$.
   \end{itemize}
   \item[(SE):] The residue field and the value group are stably embedded.
 \end{enumerate}
 
 These conditions appear implicitely in \cite[Thm.~4.1]{NTP2-00} and the previous lemmas, without being named. They are named by Jahnke and Simon in \cite{JS} in the NIP context, that is, replacing NTP2 by NIP in the (Im) condition.

 We say that a valued field has NTP2 CHIPS if it satisfies the aforementionned Chernikov-Hils' (Im) Plus (SE) conditions. Similarly, we say that a valued field has NIP CHIPS if it satisfies Jahnke-Simon's (Im) (with NIP formulas) Plus (SE). These conditions are sufficient to obtain transfer:

\begin{thm}[NTP2 CHIPS transfer, {\cite[Thm.~4.1]{CH}}]
 Let a valued field (possibly with augmented structure) $(K,v,\cdots)$ have NTP2 CHIPS, then $(K,v,\cdots)$ is NTP2 iff $(k,\cdots)$ and $(\Gamma,\cdots)$ are NTP2. 
\end{thm}

\begin{thm}[NIP CHIPS transfer, {\cite[Thm.~2.3]{JS}}] Let a valued field (possibly with augmented structure) $(K,v,\cdots)$ have NIP CHIPS, then $(K,v,\cdots)$ is NIP iff $(k,\cdots)$ and $(\Gamma,\cdots)$ are NIP. 
\end{thm}

Here we keep the assumption on $\Gamma$, since it could have augmented structure, and therefore is not automatically NIP.

By using NIP CHIPS transfer, Anscombe and Jahnke proved their classification of NIP henselian valued fields, generalizing all previously known transfer theorems.

Our goal is now to prove a \NIPn CHIPS transfer theorem and use it to achieve a classification of \NIPn henselian valued fields.

  To do so, at first in \Cref{sec-nipndefs} we recall usual definitions, theorems and conjectures around \NIPn fields. Then, in \Cref{sec-transfer}, we introduce \NIPn CHIPS and prove that it implies \NIPn transfer. The proof follows the heuristic given thereafter. With the help of this transfer theorem, we finally prove \Cref{NIPthm} in \Cref{sec-ajltr} and discuss some consequences.
  
  \subsection{Heuristic}

We give a heuristic about why having CHIPS, of any flavor, is sufficient to obtain transfer. Say the following is true for a given valued field $(K,v,\cdots)$ -- possibly with augmented structure:

\begin{enumerate}
  \item[(Im):] The type of immediate extensions is controlled by formulas of some model theoretic flavor (NIP, NTP2, etc.);
   \item[(SE):] The residue field and the value group are stably embedded.
 \end{enumerate}
 
Now assume that the (induced structure on the) residue field and value group are NIP/\NIPn/NTP2, or any other flavour, but that the valued field is not. Most of the time, combinatorial complexity can be witnessed by indiscernibles, so if a formula $\phi$ has TP2, IP or \IPn, there's a (potentially generalized) indiscernible sequence $(a_i)_{i\in I}$ and a singleton $b$ such that $\phi(b,a_i)$ witnesses some pattern. By Ramsey and compactness, we can extend $(a_i)_{i\in I}$ until each $a_i$ is a small model $K_i$. Now, of course, $\phi$ is in the type of $b$ over some $K_i$, say $K_0$ (otherwise it's always false, and that's not a pattern), but $K_0(b)/K_0$ might not be immediate. Well, whatever; let's append an enumeration of the value groups and residue field of $K_0(b)$ to $K_0$. We would like to also be able to append to the rest of the sequence $K_i$ so that the now augmented sequence keeps the indiscernability properties it had before; because the value group and residue field are stably embedded, this can be done via an array extension lemma. In the end, we have indeed that $K_0(b)/K_0$ is immediate, so $\phi$ is implied by NTP2, NIP or \NIPn formulas, and thus is itself NTP2, NIP or \NIPn -- which contradicts the choice of $\phi$.

  \subsection{Notations}\label{subsec-not}
  
  Given a valued field $(K,v)$, we write $\Gamma_v$ for its value group, $k_v$ for its residue field, $\rond O_v$ for its valuation ring and $\rond M_v$ for its maximal ideal. When the context is clear, we omit the subscript $_v$. When we consider $(K,v)$ as a first-order structure, we consider it as a 3-sorted structure, with sorts $K$ and $k$ equipped with the ring language, $\Gamma$ equipped with the ordered group language, and (partial) functions between sorts $v\colon K\rightarrow \Gamma$ and $\overline\cdot\colon K\rightarrow k$.
  
  We let lowercase letters $x$, $y$, $z$... denote variables or tuples of variables and $a$, $b$, $c$... denote parameters or tuples of parameters. We almost never use the overline to denote tuples since we prefer to let $\overline x$ be the residue of $x$ in a given valued field.
  
  Given a valued field $(K,v)$ of mixed characteristic, the standard decomposition around $p$ is defined by fixing two convex subgroups:
 
 $$\Delta_0=\bigcap_{\mathclap{\substack{v(p)\in\Delta\\\Delta\subset\Gamma\text{ convex}}}}\Delta\phantom{MM}\&\phantom{MM}\Delta_p=\bigcup_{\mathclap{\substack{v(p)\notin\Delta\\\Delta\subset\Gamma\text{ convex}}}}\Delta$$
 
 And considering the associated valuations $v_0$ and $v_p$. We then perform the following decomposition, written in terms of places with specified value groups:
 $$K\xrightarrow{\Gamma_v/\Delta_0}k_0\xrightarrow{\Delta_0/\Delta_p}k_p\xrightarrow{\Delta_p}k_v$$

We immediately remark that $\Delta_0/\Delta_p$ is of rank 1 and that $\ch(k_0)=0$ and $\ch(k_p)=p$.

We call a valued field maximal if it does not admit any proper immediate extension. Similarly, we call a valued field algebraically maximal, or separably algebraically maximal, if it does not admit any proper algebraic or proper separable algebraic immediate extension.

We call a valued field of residue characteristic $p>0$ Kaplansky if its value group is $p$-divisible and its residue field is Artin-Schreier-closed and perfect. We call all valued fields of residue characteristic 0 Kaplansky for convenience.

We shorten (separably) algebraically maximal Kaplansky in (S)AMK.

We call a mixed characteristic valued field $(K,v)$ unramified if the interval $(0,v(p))$ is empty and finitely ramified if $(0,v(p))$ contains finitely many elements.
  
We write ``we'' for ``I'', except when we write ``I'' for ``I''. I contain multiltudes.

\section{Everything you need to know about \texorpdfstring{\NIPn}{NIPn} theories}\label{sec-nipndefs}
\subsection{The \texorpdfstring{$n$}{n}-independence property}

\NIPn theories are the most natural generalization of NIP. They were first defined and studied for $n=2$ by Shelah in \cite{Shelah-strdep}. Their behavior is erratic, sometimes very similar to NIP theories, sometimes wildly different.

\begin{dfn}
 Let $T$ be a complete theory and $\rond M\satisf T$ a monster model. A formula $\phi(x;y_1,\dots,y_n)$ is said to have the independence property of order $n$ (\IPn) if there are $(a^k_i)^{1\leqslant k\leqslant n}_{i<\omega}$ and $(b_J)_{J\subset\omega^n}$ in $\rond M$ such that $\rond M\satisf\phi(b_J,a^1_{i_1},\dots,a^n_{i_n})$ iff $(i_1,\dots,i_n)\in J$.
 A formula is said to be \NIPn if it doesn't have \IPn, and a theory is called \NIPn if all formulas are \NIPn. We also write ``strictly \NIPn'' for ``\NIPn and \IP{n-1}''.
\end{dfn}

Note that having \IP{n+1} implies having \IP{n} and that \IP1 is the usual definition of the independence property (IP).

Structures which are strictly \IPn exist for all $n$:

\begin{ex}[{\cite[ex.~2.2.(2)]{CPT}}]
 The random graph is strictly \NIP2. The random $n$-hypergraph, which is the Fraïssé limit of the class of all finite $n$-hypergraphs -- which are sets of vertices equipped with a symmetrical irreflexive $n$-ary relation --, is strictly \NIPn.
\end{ex}

As for NIP, the study of \NIPn formulas can be reduced significantely by considering only atomic formulas with one singleton variable, and can also be reformulated in terms of indiscernibles -- though we only quote that result in \Cref{NIPnindi}.

\begin{prop}[{\cite[Prop.~6.5]{CPT}}]
 Being \NIPn is preserved under boolean combinations: if $\phi(x;y_1,\dots,y_n)$ and $\psi(x;y_1,\dots,y_n)$ are \NIPn, so are $\phi\et\psi$ and $\neg\phi$.
 Moreover, a theory is \NIPn iff all formulas $\phi(x,y_1,\tpp,y_n)$ with $x$ a singleton are \NIPn.\footnote{In fact, one can reduce further, and only consider the formulas with all but one variable being singletons, see \cite[Thm.~2.12]{ndep2}; however, we only require the weaker version in this article.}
\end{prop}
\subsection{Conjectures on \texorpdfstring{\NIPn}{NIPn} fields}

Starting with the celebrated work of Hempel \cite{Hem-ndep}, many known results about NIP fields have been generalized as such to \NIPn fields. In fact, the following conjecture arises naturally from work of Hempel, Chernikov, and others:

\begin{con}[The \NIPn Fields Conjecture]\label{NIPncon}
 For $n\geqslant2$, strictly \NIPn pure fields do not exist; that is, a pure field is \NIPn iff it is NIP.
\end{con}

This is for pure fields. Augmenting fields with structure -- for example by adding a relation for a random hypergraph -- will of course break this conjecture, however, ``natural'' extensions of field structure such as valuation or distinguished automorphism are believed to preserve it. In a previous article, we studied the following conjecture:
 
\begin{con}\label{NIPnhvcon}
 Strictly \NIPn henselian valued fields do not exist.
\end{con}

It is clear that \Cref{NIPnhvcon} implies \Cref{NIPncon} since the trivial valuation is henselian, and \cite[Cor.~3.14]{itme-asext} shows that in fact, these conjectures are equivalent.

We now state Shelah's conjecture for \NIPn fields:
\begin{con}[Shelah's \NIPn Conjecture]\label{SCn}
 \NIPn fields are finite, separably closed, real closed, or admit a non-trivial henselian valuation.
\end{con}

This is a famous conjecture for $n=1$, attributed to Shelah though he never wrote it down. Many results make it plausible, for example Johnson's theorem in \cite{dpfin}. If we believe in Shelah's Conjecture for $n=1$ and in the \NIPn Fields Conjecture, then we believe in Shelah's Conjecture for all $n$; however, we suggest to consider it the other way around: combining Shelah's Conjecture with the Henselian Expansion Conjecture (it's the last one I swear), we obtain a proof of the \NIPn Fields Conjecture:

\begin{con}[Jahnke's \NIPn Henselian Expansion Conjecture]\label{NIPnhensexpcon}
 Let $K$ be \NIPn as a pure field and let $v$ be a henselian valuation on $K$. Then $(K,v)$ is \NIPn as a valued field.
\end{con}

For $n=1$, this is a theorem of Jahnke, see \cite{Jah-hensexp}, which is why we name this conjecture after her. For arbitrary $n$, it is still unknown in general, but we prove it for residue characteristic $p$ in \Cref{KNIPn->KvNIPn}.

\begin{prop}
 If Shelah's \NIPn Conjecture \ref{SCn} and Jahnke's \NIPn  Henselian Expansion Conjecture \ref{NIPnhensexpcon} hold, then the \NIPn Fields Conjecture \ref{NIPncon} hold. 
\end{prop}

\begin{proof}
 Let $K$ be a \NIPn field. If it is finite, separably closed or real-closed, it is NIP. If it is neither, then by Shelah's \NIPn Conjecture, it admits a non-trivial henselian valuation. Let $v_K$ be its canonical henselian valuation. By Jahnke's Conjecture, $(K,v_K)$ is \NIPn, so the residue field $k_{v_K}$ is \NIPn. Applying Shelah's Conjecture to $k_{v_K}$, it is in turn finite, separably-closed, real-closed, or admit a non-trivial henselian valuation. But by definition of $v_K$, $k_{v_K}$ can only admit non-trivial henselian valuations if it is separably closed. So, it is either separably closed, real-closed, or finite; in all cases, it is NIP, and \cite[Cor.~3.13]{itme-asext} implies that $(K,v_K)$ is NIP -- and thus $K$ is NIP.
\end{proof}

This gives a strategy for proving the \NIPn Fields Conjecture: based on the case $n=1$, prove Jahnke's Henselian Expansion Conjecture for arbitrary $n$, then wait for a proof of Shelah's Conjecture for NIP fields, and generalizes that proof to the \NIPn context.

\section{\texorpdfstring{\NIPn}{NIPn} transfer}\label{sec-transfer}

\subsection{\texorpdfstring{\NIPn}{NIPn} \& generalized indiscernibles}
\begin{dfn}
  Let $\rond{M}$ be an $\rond{L}$-structure and $\rond{I}$ be an $\rond{L}_0$-structure, where $\rond{L}$ and $\rond{L}_0$ are possibly different languages. A sequence $(a_i)_{i \in I}$ of tuples of $M$ is said to be $\rond{I}$-indiscernible over a set $A\subset\rond M$ if for any $i_0,\tpp,i_n$ and $j_0,\tpp,j_n$ in $\rond I$, $\qftp_{\rond{L}_0}(i_0,\tpp,i_n)=\qftp_{\rond{L}_0}(j_0,\tpp,j_n)$ implies $\tp_{\rond{L}}(a_{i_0},\tpp,a_{i_n}/A)=\tp_{\rond{L}}(a_{j_0},\tpp,a_{j_n}/A)$.
\end{dfn}

\begin{rem}
 We call two tuples of elements of a structure $a$ and $b$ ``of the same mould'' if they are of the same length $n$ and if for all $i<n$, $a_i$ and $b_i$ are in the same sort. Given a tuple $a$, we say that a tuple of variable $x$ is ``a mould'' of $a$ if they are of the same length $n$ and for all $i<n$, $x_i$ is a variable on the sort containing $a_i$. A contrario, given a tuple of variables $x$, we say that a tuple of elements $a$ is ``a cast'' of $x$ if $x$ is a mould of $a$, and we say similarily that two tuples of variables $x$ and $y$ are ``identical as moulds'' if $x$ is the mould of a cast of $y$.
 
 The raison d'être of these notions is to make clear that there's no reason an arbitrary sequence $(a_i)_{i\in I}$ has to be a sequence of tuples of the same mould. For a generalized indiscernible sequence, we do not need to compare the type of $a_i$ and $a_j$ if $i$ and $j$ have different types, so they might as well be of different lengths and of different sorts. This notably happens when we work with sequences indexed by partitioned structures.
\end{rem}

We denote by $G_n$ a countable ordered random $n$-partite $n$-hypergraph; it is a structure in the language $\set{<,P_1,\tpp,P_n,R}$, where $<$ is a binary relation symbol, $P_i$ are unary predicates, and $R$ is an $n$-ary relation symbol, and its complete theory is axiomatized as follows:
\begin{enumerate}
 \item $G_n=P_1\sqcup\cdots\sqcup P_n$,
 \item $<$ is a dense linear order without endpoints on each $P_i$,
 \item $P_1<\cdots<P_n$,
 \item $R$ is an $n$-ary relation on $P_1\times\tpp\times P_n$ -- the ``hyperedge'' relation,
 \item For any finite disjoint $A_0,A_1\subset P_1\times\tpp\times P_{j-1}\times P_{j+1}\times\tpp\times P_n$ and for any $b_0<b_1\in P_j$, there is $b\in P_j$ such that $b_0<b<b_1$ and for any $(g_1,\tpp,g_{j-1},g_{j+1},\tpp,g_n)\in A_0$, then $(g_1,\tpp,g_{j-1},b,g_{j+1},\tpp,g_n)$ is an edge; and same goes for $A_1$ with non-edges.
\end{enumerate}

We say that $(g_1,\tpp,g_n)$ is an edge to signify $G_n\satisf R(g_1,\tpp,g_n)$. In particular, it implies $g_i\in P_i$.

We denote by $O_n$ the reduct of $G_n$ to the language $\set{<,P_1,\tpp,P_n}$, that is, we ignore the edges. The complete theory of $O_n$ is axiomatized by axioms 1 to 3 above.

\begin{figure}[!htbp]
\begin{center}
\begin{tikzpicture}
\foreach \i in {1,2,3,4}{ 
\draw[very thick] (2*\i,0)--(2*\i,6);
\node(P\i) at (2*\i,6.4){$P_{\i}$};
}
\node[circle,fill=black,inner sep=1.5,label={\ \ \ \ $b_0$}](b0) at (6,2) {};
\node[circle,fill=black,inner sep=1.5,label={\ \ \ \ $b_1$}](b1) at (6,5) {};
\node[circle,fill=black,inner sep=1.5,label={\ \ \ \ $b$}](b) at (6,4) {};
\node[circle,fill=red,inner sep=1.5](a11) at (2,2.5){};
\node[circle,fill=red,inner sep=1.5](a12) at (2,3.1){};
\node[circle,fill=red,inner sep=1.5](a13) at (2,5){};

\node[circle,fill=green!50!black,inner sep=1.5]() at (2,2){};

\node[circle,fill=red,inner sep=1.5](a21) at (4,3){};
\node[circle,fill=red,inner sep=1.5](a22) at (4,5.2){};

\node[circle,fill=green!50!black,inner sep=1.5]() at (4,1.5){};
\node[circle,fill=green!50!black,inner sep=1.5]() at (4,3.5){};
\node[circle,fill=green!50!black,inner sep=1.5]() at (4,4){};

\node[circle,fill=red,inner sep=1.5](a31) at (8,4.6){};

\node[circle,fill=green!50!black,inner sep=1.5]() at (8,3.3){};
\node[circle,fill=green!50!black,inner sep=1.5]() at (8,2.3){};

\foreach\i in{1,2,3}{
\foreach\j in{1,2}{
\draw[red] (a1\i)--(a2\j)--(b)--(a31);
}}
\end{tikzpicture}
\end{center}
\caption[An ordered random 4-hypergraph.]{An ordered random 4-hypergraph. Each $P_i$ is represented by a vertical line. Sets $A_0$ and $A_1$ are represented in red and in green respectively. Edges are drawn in red. Such a graph will have many more edges which are not drawn here, and $A_0$, $A_1$ need not be product sets in general.}
\end{figure}
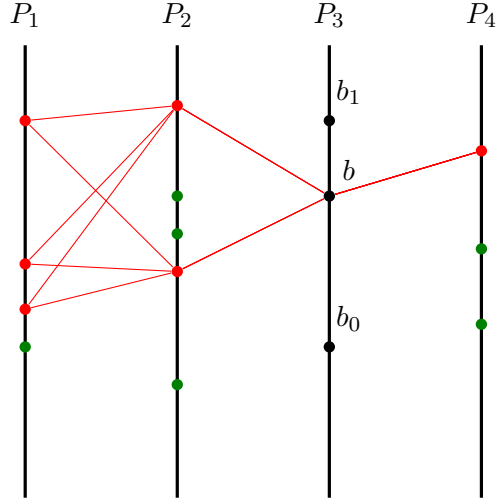

\begin{prop}[{\cite[Prop.~5.2]{CPT}, \cite[Prop.~2.8]{ndep2}}]\label{NIPnindi}
  A formula $\phi(x;y_1,\tpp,y_n)$ has \IPn iff there exists (in a sufficiently saturated model $\rond{M}$) a tuple $b$ and a sequence $(a_g)_{g\in G_n}$ which is $O_n$-indiscernible over $\emptyset$ and $G_n$-indiscernible over $b$ such that $\phi(b;y)$ encodes the edges of the graph; that is:
  $$\rond{M}\satisf\phi(b,a_{g_1},\tpp,a_{g_n})\text{ iff }G_n\satisf R(g_1,\tpp,g_n).$$
\end{prop}

Note that considering a sequence indexed by $G_n$ which is $O_n$-indiscernible is the same as considering $n$ mutually indiscernible sequences indexed by each $P_i$.

\subsection{\texorpdfstring{\NIPn}{NIPn} CHIPS transfer}

We now prove that \NIPn CHIPS implies \NIPn transfer. To do so, we first need to obtain an array extension lemma. We do so in an arbitrary complete theory $T$ with a given monster model $\rond M$.

\begin{dfn}
 Let $D$ be a $\emptyset$-definable set. We say that $D$ is $n$-hanced stably embedded if for all formulas $\phi(x,y_1,\tpp,y_n)$ and for all sequences $(a_i^k)_{i\in I}^{1\leqslant k\leqslant n}\in\rond M$ such that each $a_i^k$ is a cast of $y_k$, there is a formula $\psi(x,z_1,\cdots,z_n)$ and a sequence $(b_i^k)_{i\in I}^{1\leqslant k\leqslant n}\in D$ -- with each $b_i^k$ a cast of $z_k$ -- such that:
 $$\phi(D,a_{i_1}^1,\tpp,a_{i_n}^n)=\psi(D,b_{i_1}^1,\tpp,b_{i_n}^n).$$
\end{dfn}

The usual definition of stable embeddedness is that any $\rond M$-definable subset of $D$ is $D$-definable. A priori, this $D$-definition depends wildly on the original $\rond M$-definition, however, with compactness and coding tricks, this can be strengthened to a uniform version. This is discussed in great detail in \cite[sec.~1]{SE-pierre}.

Our version is semi-uniform -- $\psi$ depends on $\phi$ and also on the choice of the sequence $(a_i)_{i\in I}$, but does not change when going from $a_i$ to $a_{i'}$ --, and more importantly, it works on $n$ variables at once. It might be that this is equivalent to being stably embedded, assuming $D$ is infinite, via a coding trick and a compactness argument; but it remains to be proved. We note the following:

\begin{lem}\label{autoSE}
 If every automorphism of $D^n$ lifts to an automorphism of $\rond M^n$, then $D$ is $n$-hanced stably embedded.
\end{lem}

\begin{proof}
 This can be obtained by adapting the proof of \cite[App.~Lem.~1]{CH-SE}, specifically, the proof of (6) implies (5). Note that if $D$ is not $n$-hanced stably embedded, then there exists an $\rond M$-definable family $S_{a_1,\tpp,a_n}=\set{b\in D}{\rond M\satisf\phi(b,a_1,\tpp,a_n)}$ which is not a $D$-definable family. Following the original proof with this definable family instead of a mere definable set yields the wanted result.
\end{proof}

In order to study $n$-hanced stable embeddedness in more detail, we ideally would want an $n$-hanced version of the aforementioned lemma \cite[App.~Lem.~1]{CH-SE}, this has not been achieved as of yet.

\begin{lem}
  Let $(a_g)_{g\in G_n}$ be $O_n$-indiscernible over a set $A$. Suppose $D$ is a $\emptyset$-definable set which is $n$-hanced stably embedded and fix $d\in D$. If the induced structure on $D$ is \NIPn, then no formula with parameters in $Ad$ can encode the edges of $(a_g)_{g\in G_n}$.
\end{lem}

This is a \NIPn version of \cite[Lem.~2.1]{JS}.

\begin{proof}
  Let $\phi(d;y_1,\!\cdots\!,y_n)$ be a formula with unwritten parameters in $A$ and encoding the edges of $(a_g)_{g\in G_n}$. By $n$-hanced stable embeddedness, we can find $\psi(x,z_1,\!\cdots\!,z_n)$ and $(b_g)_{g\in G_n}\in D$ such that $\phi(D;a_{g_1},\!\cdots\!,a_{g_n})=\psi(D;b_{g_1},\!\cdots\!,b_{g_n})$ for all $(g_1,\tpp, g_n)$. 
  
  \begin{claim} For any $J\subset P_1\times\cdots\times P_n$, we can find $d_J\in D$ such that $\phi(d_J;a_{g_1},\!\cdots\!,a_{g_n})$ holds iff $(g_1,\!\cdots\!,g_n)\in J$.
  \end{claim}

Given such $d_J$, we immediately have that $\psi(d_J;b_{g_1},\!\cdots\!,b_{g_n})$ holds iff $(g_1,\!\cdots\!,g_n)\in J$, which yields \IPn on $D$; thus proving the claim is enough to prove the lemma.

To prove the claim, let $f$ enumerate $P_1\times\cdots\times P_n$ in such a way that $f(i)$ and $f(i+1)$ always differ in exactly one coordinate -- remember that $G_n$ is countable. We will prove that one can find a $d_N\in D$ such that $\phi(d_N;a_{f(i)_1},\tpp,a_{f(i)_n})$ holds iff $f(i)\in J$ for $i<N$. For $N=1$, either $f(0)$ is in $J$ or not. We can find $(g_1,\!\cdots\!,g_n)\in P_1\times\cdots\times P_n$ such that $\phi(d;a_{g_1},\!\cdots\!,a_{g_n})$ holds (or not), so $M\satisf\exists x\in D(\neg)\phi(d;a_{g_1},\!\cdots\!,a_{g_n})$, and by $O_n$-indiscernability, $M\satisf\exists x\in D(\neg)\phi(d;a_{f(0)_1},\tpp,a_{f(0)_n})$.

Now assume such a $d_N$ exists for some $N$. We do the case $f(N)\in J$, the other one is similar. We need to find $(g_1,\!\cdots\!,g_n)$ in the same place as $f(N)$ regarding $f(i),i<N$ (that is, $g_j<f(i)_j$ iff $f(N)_j<f(i)_j$, etc.), forming an edge, and not colliding with previous choices.

By our choice of $f$, $f(N-1)_i=f(N)_i$ for all $i$ but 1. Then we take $g_i=f(N-1)_i$, and for the remaining $g_j$, we use the properties of $G_n$:

\ 

If $f(N)_j$ has never appeared before, we just need to choose a $g_j$ in the correct place such that $(g_1,\!\cdots\!,g_n)$ forms an edge. This is possible by the properties of the random hypergraph.

\ 

If $f(N)_j$ has appeared before, then fixing $g_j=f(N)_j$ might cause trouble, since $(g_1,\!\cdots\!,g_n)$ might not be connected. Instead, we let $I$ be the set of $i$ such that $f(i)_j=f(N)_j$. We define $b_0$ and $b_1$ in $P_j$ as follows:
$$b_0=\max\set{f(i)_j}{f(i)_j<f(N)_j}\text{ and }b_1=\min\set{f(i)_j}{f(i)_j>f(N)_j}\text{.}$$ We let $A_0$ be the set of $(f(i)_1,\!\cdots\!,f(i)_{j-1},f(i)_{j+1},\!\cdots\!,f(i)_n)$, for $i\leqslant N$, such that $$(f(i)_1,\!\cdots\!,f(i)_{j-1},f(N)_j,f(i)_{j+1},\!\cdots\!,f(i)_n)$$ forms an edge. We let $A_1$ be the counterpart with non-edges. Then by the properties of $G_n$, there is $b$ between $b_0$ and $b_1$, forming edges with all points of $A_0$ and no points of $A_1$; we now let $f'(i)=f(i)$ for $i\notin I$, $f'(i)_k=f(i)_k$ for $k\neq j$, and $f'(i)_j=b$ for $i\in I$. We conclude by indiscernability as before.
 
 This proves the claim, and thus the lemma.
\end{proof}

\begin{lem}[\NIPn array extension lemma]\label{interpol}
  Let $D$ be $n$-hanced stably embedded and let $(a_g)_{g\in G_n}$ be $O_n$-indiscernible over $\emptyset$ and $G_n$-indiscernible over some tuple $b$. Fix an edge $(g_1,\tpp,g_n)\in P_1\times\tpp\times P_n$. For each $g_i$ let $c_{g_i}\in D$ be a small tuple. Then, we can interpolate the rest of the sequence, that is, we can find $(c_g)_{g\neq g_i}$ and $(a'_g)_{g\in G}$ such that:
  \begin{itemize}
    \item $a'_{g_i}=a_{g_i}$,
    \item $\tp((a'_g)_{g\in G_n}/b)=\tp((a_g)_{g\in G_n}/b)$,
    \item $(a'_g c_g)_{g\in G_n}$ is $O_n$-indiscernible over $\emptyset$ and $G_n$-indiscernible over $b$.
  \end{itemize}
\end{lem}

This is a \NIPn version of \cite[Lem.~2.2]{JS} and \cite[Lem.~3.8]{CH}.

\begin{proof}
    We do it part by part, mimicking the strategy of the NIP case. We fix an edge $(g_1,\tpp,g_n)\in G_n$, and we fix $i$. In the NIP case, we do even and odd separately; here we define the set of ``even'' indices to be $$E_i=\set{g\in P_i}{(g_1,\tpp g_{i-1},g,g_{i+1},\tpp,g_n)\text{ is an edge}}.$$ Because $(a_g)_{g\in G_n}$ is $G_n$-indiscernible over $b$, we can find $c_g$ for each $g\in E_i$ such that $$a_gc_g\equiv_{b,a_{g_1},\tpp,a_{g_{i-1}},a_{g_{i+1}},\tpp,a_{g_n}}a_{g_1}c_{g_1}.$$ Now, by Ramsey, we may assume $(a_gc_g)_{g\in G_n,g\notin(P_i\setminus E_i)}$ is $O_n$-indiscernible over $\emptyset$ and $G_n$-indiscernible over $b$.
    
    Now, because this is true for any sequence with these properties, we move to a new sequence where $P_i$ is now $P_i^*$ and is very long. Any ``even'' element of $P_i^*$ has already been extended by a $c$.
    
    For each element of $g\in E_i$ (the original, short version) we chose a representation $\lambda_g\in P_i^*$. We make sure to take them very far apart from each other.
    
    $P_i\setminus E_i$ injects into the set of cuts of $E_i$. Fix an ``odd'' index $h$, and look at the corresponding cut $C_h$ (in $P^*_i$) of $\set{\lambda_g}{g\in E_i}$. Now $P_1\sqcup\tpp\sqcup P_{i-1}\sqcup C_h\sqcup P_{i+1}\sqcup\tpp\sqcup P_n$ is itself a random graph.
    
    Take a formula $\phi(a_{g_i},c_{g_i})\in\tp(a_{g_i}c_{g_i}/b(a_g)_{g\notin P_i})$. By the previous lemma, since $\phi$ can't encode the random graph, $\phi(a_k,c_k)_{k\in C_h}$ must be either always true or always false, except for discretely many points.
    
    If we exclude those discretely many points from $C_h$, after having done that for all formulas, we still have points, because $P_i^*$ is really long. Chose any ``even'' point in what remains; we will call it $\lambda_h$.
    
    Now we take an automorphism $\sigma$ over $b(a_g)_{g\notin P_i}$ taking each $a_{\lambda_g}c_{\lambda_g}$ to $a_gc_g$. We define $a'_hc_h=\sigma(a_{\lambda_h}c_{\lambda_h})$. Now the sequence with extended points in the $i$th part and $a'$ for ``odd'' indices satisfy the theorem.
\end{proof}

 We now suppose $T$ is a complete theory of valued fields (possibly with additional structure), and we consider the following properties:
 
 \begin{enumerate}
   \item[(SE)$_n$:]\label{SE} The residue field and the value group are $n$-hanced stably embedded.
   \item[(Im)$_n$:]\label{Im} For any small models $K_1,\tpp, K_n\satisf T$, writing $L$ for the compositum of all of them, and for any singleton $b\in\rond{M}$, if $L(b)/K_i$ is immediate for all $i$, then we have that $\tp(b/K_1,\tpp,K_n)$ is implied by instances of \NIPn formulas, that is, there is a $p\subset\tp(b/K_1,\tpp,K_n)$ such that:
   \begin{itemize}
    \item any formula $\phi(x;y_1,\tpp,y_n)\in p$ -- where $x$ is the cast for $b$ and $y_i$ for $K_i$ -- is \NIPn, and
    \item $\psi(b,K_1,\tpp,K_n)$ holds iff $p\vdash\psi$.
   \end{itemize}
 \end{enumerate}
 
 We say that (the complete theory of) a valued field, potentially with augmented structure, has \NIPn CHIPS if it checks these two conditions.

 \begin{cor}[\NIPn CHIPS transfer]\label{SEIMtrns}\label{NIPnchips}
  If $T$ is a complete theory of valued fields with \NIPn CHIPS, then $T$ has \NIPn transfer; that is, $T$ is \NIPn iff the theories of the residue field and the value group are \NIPn.
\end{cor}

This is a \NIPn version of \cite[2.3]{JS}. Let us also note that in the case where the structure is augmented, when checking whether a theory has CHIPS -- whether it be of NIP, \NIPn or NTP2 flavour --, we need to be careful on exactly what is the structure we consider on the residue field and on the value group; if for example $k$ is \NIPn as a pure field, but we only know that an augmented structure of $k$ is (SE)$_n$, augmented structure for which we don't know \NIPn, then this theorem does not guarantee transfer.

\begin{proof}
  Assume $T$ has \IPn. Then we can find a formula $\phi(x;y_1,\tpp,y_n)$ with $x$ unary, a singleton $b$ and a sequence $(a_g)_{g\in G_n}$ $O_n$-indiscernible over $\emptyset$ and $G_n$-indiscernible over $b$, such that $\phi(b;a_{g_1},\tpp,a_{g_n})$ holds iff $G_n\satisf R(g_1,\tpp,g_n)$.
  
  By Ramsey and compactness, we can extend each $a_g$ until it enumerates a small model $K_g$. We refer to \cite{CPT}, specifically the appendix, for the study of Ramsey properties in \NIPn theories.
  
  We fix an edge $(g_1,\tpp,g_n)\in P_1\times\tpp\times P_n$. Let $k'$ and $\Gamma'$ be the residue and value group of $K_{g_1}\cdots K_{g_n}(b)$, let $c_{g_i}$ and $d_{g_i}$ be enumerations of $k'\setminus k_{g_i}$ and $\Gamma'\setminus\Gamma_{g_i}$. Apply the previous lemma twice to obtain a sequence $(a'_g c_g d_g)_{g\in {G_n}}$ such that:
  \begin{itemize}
    \item $a'_{g_i}=a_{g_i}$,
    \item $\tp((a'_g)_{g\in G_n}/b)=\tp((a_g)_{g\in G_n}/b)$,
    \item $(a'_g c_g d_g)_{g\in G_n}$ is $O_n$-indiscernible over $\emptyset$ and $G_n$-indiscernible over $b$.
  \end{itemize}
  We now start over: we extend each $(a'_g c_g d_g)$ to enumerate a model, add the residue and value group of this model plus $b$, and interpolate.
  After $\omega$ iterations, we have a sequence $(N_g)$ of small models, $O_n$-indiscernible over $\emptyset$, $G_n$-indiscernible over $b$, such that $\tp((N_g)_{g\in G_n}/b)$, restricted to the correct subtuple, equals $\tp((a_g)_{g\in G_n}/b)$, and such that $N_{g_1}\cdots N_{g_n}(b)/N_{g_i}$ is immediate.
  Now by (Im)$_n$, $\tp(b/N_{g_1},\tpp,N_{g_n})$ is implied by instances of \NIPn formulas. By $G_n$-indiscernability, such a formula will also hold for any edge. But by \NIPn-ity, it can't also not hold for all non-edges, in fact it can only not hold for finitely many of them.
  Hence we must have a non-edge $(g'_1,\tpp,g'_n)$ such that all the \NIPn formulas implying $\tp(b/N_{g_1},\tpp,N_{g_n})$ hold, and thus $\phi(b,a_{g'_1},\tpp,a_{g'_n})$ holds, which contradicts the initial choices of $\phi$, $b$, and $a$.
\end{proof}

\section{A \texorpdfstring{\NIPn}{NIPn} Anscombe-Jahnke}\label{sec-ajltr}

We now use \Cref{SEIMtrns} to prove \Cref{NIPthm}.

\begin{prop}\label{SAMKSEn}
 SAMK henselian valued fields have (SE)$_n$.
\end{prop}

\begin{proof}
 By \Cref{autoSE}, it is enough to show that every automorphism of $\Gamma^n$ lifts to $K^n$, and similarly for every automorphism of $k^n$. This follows directly from adapting the proof of Anscombe-Jahnke in the case $n=1$, see \cite[Thm.~12.6]{AJ-Cohen}.
\end{proof}

\begin{prop}\label{SAMKtrans}
 If $(K,v)$ is SAMK with \NIPn residue, then it is \NIPn.
\end{prop}

\begin{proof}
 For $n=1$, this was done by Jahnke and Simon in the case of finite degree of imperfection, and Anscombe and Jahnke for the rest; see \cite[Thm.~3.3]{JS} and \cite[Prop.~4.1]{AJ-NIP}.
 
 The previous proposition tells us $(K,v)$ has (SE)$_n$, we now prove it has (Im)$_n$: let $K_1,\tpp, K_n$ be small models of the theory of $(K,v)$ -- as always we are working in a monster model, thus all valuations are restriction of a given valuation on the monster -- and $b$ a singleton such that $K_1\tpp K_n(b)/K_i$ is immediate. We let $L$ be the henselization of the relative perfect hull of $K_1\tpp K_n(b)$. By the properties of the henselization, $L$ is uniquely determined by the isomorphism type of $b$ over $K_1\tpp K_n$.
 
 Now we consider $L'$, the relative tame closure of $L$. This is uniquely determined up to isomorphism by \cite[Thm.~5.1]{KPR} because $L$ is Kaplansky. By \cite[Thm.~5.1]{Delths}, $L'$ is an elementary extension of $K_i$ (for any $i$).
 
 Thus, the isomorphism type of $b$ over $K_1, \tpp, K_n$ (that is, its qf type) uniquely determines a model containing it, so it implies the full type. Quantifier free formulas in the language of valued fields are NIP, thus in particular \NIPn; which means $(K,v)$ has (Im)$_n$, and we have transfer by \Cref{SEIMtrns}.
\end{proof}

Note that we did not specify the characteristic -- the way we wrote it assumes the residue characteristic is $p$, but in equicharacteristic 0, it's even simpler, since $K_1\tpp K_n(b)\equiv k((\Gamma))\equiv K_i$.

In equicharacteristic, we already proved that \NIPn henselian valued fields are SAMK (or trivial), so this suffices to have the equivalence, and only the mixed characteristic case remains.

\begin{lem}
 If $(K,v)$ is henselian, of mixed characteristic and unramified, then it has (SE)$_n$.
\end{lem}

\begin{proof}
 As before, it is an easy adaptation of the proof in the case $n=1$, see \cite[Lem.~3.1]{JS} and \cite[Prop.~4.1]{AJ-NIP}, using \Cref{autoSE}.
\end{proof}

\begin{lem}\label{unramtr}
 If $(K,v)$ is (mixed-char) unramified with \NIPn residue, then it is \NIPn.
\end{lem}

\begin{proof}
 Again, \NIP1-transfer has been proved using (SE)$_1$+(Im)$_1$ by Anscombe and Jahnke, see \cite[Lem.~4.4]{AJ-NIP}. We now go towards arbitrary $n$.
 
 We let $K_1,\tpp,K_n$ be small models -- of a given monster model, as above -- and $b$ be a singleton such that $K_1\tpp K_n(b)/K_i$ is immediate for each $i$. We also assume that one of them, say $K_1$, is $\aleph_1$-saturated. Each of them is equipped with a valuation which is the restriction of the monster's valuation and that we denote $v$ in each of them.
 
 Let $L=K_1\tpp K_n(b)$, by assumption $L/K_1$ is immediate, so we write $\Gamma$ for the value group and $k$ for the residue field. By unramification, $\Gamma=\Delta\oplus\mathds Z$, with $\Delta=\Gamma/\mathds Z$ and $v(p)=(0,1)\in\Delta\oplus\mathds Z$, and we let $w$ be the coarsening of $v$ corresponding to $\mathds Z$. We denote the residue field of $(\cdot,w)$ by $\overline{\cdot}$.
 
 Now $(\overline{L},\overline{v})$ is an immediate extension of $(\overline{K_1},\overline{v})$. But by $\aleph_1$-saturation, $(\overline{K_1},\overline{v})$ is spherically complete, hence maximal. So, $\overline{L}=\overline{K_1}$.
 
 Finally, we consider the henselization $L^h$ of $L$. It is immediate over $L$ – and over $K_1$. Decomposing it into its $\Delta$ part and its $\mathds Z$ part, we have that $\overline{L^h}=\overline{L}^h=\overline{L}$, since it is equal to $\overline{K_1}$ which is henselian.
 \begin{center}
 \begin{tikzpicture}[node distance=2cm, on grid]
  \node (K) {$K_1$};
  \node (Kbar)[right=of K] {$\overline {K_1}$};
  \node (k)[right=of Kbar] {$k$};
  
  \node (L)[above=of K] {$L$};
  \node (Lbar)[right=of L] {$\overline L$};
  \node (l)[right=of Lbar] {$k$};
  
  \node (Lh)[above=of L] {$L^h$};
  \node (Lhbar)[right=of Lh] {$\overline {L^h}$};
  \node (lh)[right=of Lhbar] {$k$};
  
  \draw[->,above](K) to node {$\Delta$} (Kbar);
  \draw[->,above](L) to node {$\Delta$} (Lbar);
  \draw[->,above](Lh) to node {$\Delta$} (Lhbar);
  
  \draw[->,above](Kbar) to node {$\mathds Z$} (k);
  \draw[->,above](Lbar) to node {$\mathds Z$} (l);
  \draw[->,above](Lhbar) to node {$\mathds Z$} (lh);
  
  \draw[-](K) to node {}(L) to node {}(Lh);
  \draw[-,right](Kbar) to node {$=$}(Lbar) to node {$=$}(Lhbar);
 \end{tikzpicture}
 \end{center}
 
 The $\mathds Z$ part of $L^h$ and $K_1$ are exactly the same, this implies that $(K_1,v)$ is an elementary substructure of $(L^h,v^h)$ by \cite[Cor.~12.5]{AJ-Cohen}
 
 This means that the quantifier free type of $b$ over $K_1$ completely determines a model containing $K_1\cdots K_n(b)$, that is, it implies the full type $\tp(a/K_1,\tpp,K_n)$. Note that we fixed $K_1$ but we could have worked over any $K_i$ instead.
\end{proof}

We need to go from unramified to finitely ramified, and to study compositions of valuations in the standard decomposition. The following results will be useful:

\begin{prop}\label{SEaugment}
    Let $\rond{L}$ be relational, let $\rond M$ be a \NIPn $\rond L$-structure, let $D$ be $\emptyset$-definable and $n$-hanced stably embedded. Consider an extension $D'$ of $D_{ind}$ to a relational language $\rond L_p$, and let $\rond M'$ be the corresponding extension of $\rond M$ to $\rond L'=\rond L\cup \rond L_p$.
    
    Then, $D'$ is $n$-hanced stably embedded in $\rond M'$, and if furthermore $D'$ is \NIPn, then so is $\rond M'$.
\end{prop}

Before proving it, let us specify how we will use it: we aim to obtain a \NIPn version of \cite[Prop.~3.3]{AJ-NIP}. To do so, we apply the proposition above with $\rond{L}$ a relational version of the language of valued fields, $\rond M=(K,w)$, $D=k_w$, and $\rond L_p$ containing a predicate for a valuation $\overline v$ on $D=k_w$, and we get:

\begin{cor}\label{NIPncompo}
    Let $(K,v)$ be a valued field and $w$ be a coarsening of $v$. Assume that $(K,w)$ and $(k_w,\overline v)$ are both \NIPn and that $k_w$ is $n$-hanced stably embedded (as a pure field) in $(K,w)$. Then $(K,v)$ is \NIPn.
\end{cor}

\begin{proof}[Proof of \Cref{SEaugment}.]
    We may assume that $D'$ has QE in $\rond L_p$ and $\rond M$ in $\rond{L}$; then (the proof of) \cite[Lem.~46]{CS15} implies that every $\rond L'$-formula is equivalent to a $D$-bounded formula, that is, a formula of the form:
    $$Qy\in D\Ou_{i<m}\phi_i(x, y)\et \psi_i(x, y)$$
    with $Q$ a tuple of quantifiers, $\phi_i$ qf-$\rond L$-formulas and $\psi_i$ qf-$\rond{L}_p$-formulas (with $x$ restricted to $D$).
    
    Thus, $D'$ is $n$-hanced stably embedded in $\rond M'$, and its induced structure is exactly coming from $\rond{L}_p$.
    
    We now assume $D'$ is \NIPn and we prove by induction on the number of quantifiers that every $D$-bounded formula is \NIPn. If it has no quantifier, it is \NIPn by assumption. Now let $\phi(x,y_1,\tpp,y_n)=\exists z\!\in\!D\,\psi(x,y_1,\tpp,y_nz)$, where $\psi$ is $D$-bounded and \NIPn.
    
    Suppose $\phi$ has \IPn. Then, in a sufficiently saturated model, we can find $(a_g)_{g\in G_n}$ and $b$ such that $(a_g)_{G_n}$ is $G_n$-indiscernible over $b$ and $O_n$-indiscernible over $\emptyset$. Fix an edge $(g_1,\tpp, g_n)$, now $\exists z\in D\psi(b,a_{g_1},\tpp,a_{g_n}z)$ holds and we can find $c_{g_n}\in D$ witnessing it. Interpolate the sequence using \Cref{interpol} to get $(a'_g)_{g\in G_n}$ and $(c_g)_{g\in G_n}$ such that:
      \begin{itemize}
    \item $a'_{g_i}=a_{g_i}$,
    \item $\tp((a'_g)_{g\in G_n}/b)=\tp((a_g)_{g\in G_n}/b)$,
    \item $(a'_g c_g)_{g\in G_n}$ is $O_n$-indiscernible over $\emptyset$ and $G_n$-indiscernible over $b$.
  \end{itemize}
  By $G_n$-indiscernability over $b$, since $\psi(b,a_{g_1},\tpp,a_{g_n}c_{g_n})$ holds, it also holds for any edge. By assumption, $\forall z\in D\neg\psi(b,a_{g'_1},\tpp,a_{g'_n}z)$ holds for any non-edge, thus in particular not for $z=c_{g'_n}$.
  
  Hence there is an \IPn pattern for $\psi$, which contradicts our induction hypothesis.
\end{proof}

\begin{prop}\label{finramtrans}
 Let $(K,v)$ be henselian of mixed characteristic such that $(K,v_p)$ is finitely ramified and $(k_p,\overline v)$ is \NIPn; then $(K,v)$ is \NIPn.
\end{prop}

\begin{proof}
 Since $v_p$ is finitely ramified, it is definable; see \cite[Cor.~1.4.3]{itmePHD}. Thus, if we consider an $\aleph_1$-saturated extension $(K^*,v^*)$ of $(K,v)$, we have that $(K^*,v_p^*)$ is also finitely ramified, and $(k_p^*,\overline v^*)$ is also \NIPn. Furthermore, $(K,v)$ is \NIPn iff $(K^*,v^*)$ is \NIPn; thus we may assume that $(K,v)$ is $\aleph_1$-saturated.
 
 As usual, we consider the standard decomposition. By $\aleph_1$-saturation, $(k_0,\overline v_p)$ is complete; it is also rank-1 by definition and finitely ramified by assumption. By \cite[Thm.~22.7]{war-toporin}, there is a field $L$ such that $k_0/L$ is finite and such that, writing $w=\overline {v_p}\vert_L$, we have that $(L,w)$ is complete, unramified, and has residue field $k_w=k_p$.
 
 \begin{center}
  \begin{tikzpicture}
   \node(K) {$K$};
   \node(k0)[right=of K]{$k_0$};
   \node(kp)[right=of k0]{$k_p$};
   \node(k) [right=of kp]{$k$};
   \node(L) [below=of k0]{$L$};
   
   \draw[->,above] (K) to node {$v_0$} (k0);
   \draw[->,above] (k0) to node {$\overline{v_p}$} (kp);
   \draw[->,above] (kp) to node {$\overline v$} (k);
   \draw[->,below right] (L) to node {$w$} (kp);
   \draw[-,left] (L) to node {finite} (k0);
  \end{tikzpicture}
 \end{center}

 Since we know that $k_p$ is \NIPn, by \Cref{unramtr}, $(L,w)$ is \NIPn; we also know that $k_p$ is $n$-hanced stably embedded in $(L,w)$. We are thus in the setting of \Cref{NIPncompo}, so $(L,\overline v\circ w)$ is \NIPn. Since $k_0$ is a finite extension of $L$, we conclude that s$(k_0,\overline v)$ is \NIPn as well.
 
 Finally, we apply \Cref{unramtr} once more to the fields $(K,v_0)$ and $(k_0,\overline v)$: because $(K,v_0)$ is of equicharacteristic 0, $k_0$ is $n$-hanced stably embedded, and since it is \NIPn, we know $(K,v_0)$ is \NIPn by equicharacteristic 0 transfer. Since we just proved that $(k_0,\overline v)$ is \NIPn, $(K,v)$ itself is \NIPn.
\end{proof}

We are now finally ready to prove our main theorem.

\begin{proof}[Proof of \Cref{NIPthm}.] Let $(K,v)$ be a henselian valued field.

If $(K,v)$ is \NIPn, then so is its residue field $k$. For the rest, we conclude by \Cref{asext}. This gives one direction of the theorem.

In the other direction, assume that $k$ is \NIPn. If $v$ is trivial then $(K,v)$ is \NIPn. Assume $v$ is non-trivial. If $K$ is of equicharacteristic and SAMK, then $(K,v)$ is \NIPn by \Cref{SAMKtrans}. If $K$ is of mixed characteristic, $(K,v_p)$ finitely ramified, and $(k_p,\overline v)$ SAMK or trivial; then $(K,v)$ is \NIPn by \Cref{finramtrans}. Finally, if $K$ is of mixed characteristic and $(k_0,\overline v)$ is AMK, then $(k_0,\overline v)$ is \NIPn by \Cref{SAMKtrans} -- since AMK and SAMK are the same thing for a characteristic 0 field such as $k_0$. Finally, we conclude that $(K,v)$ is \NIPn by applying \Cref{NIPncompo}: $(K,v_0)$ is of equicharacteristic 0 so $k_0$ is $n$-hanced stably embedded in it, $(k_0,\overline v)$ is \NIPn, hence $(K,v)$ is \NIPn.
\end{proof}

\begin{cor}\label{NIPn+NIPres->NIP}
 Let $(K,v)$ be a henselian valued field. Assume $K$ is \NIPn. If $k_v$ is \NIP{m} for some $m<n$, then $(K,v)$ is \NIP{m}. In particular, if $k$ is NIP, $(K,v)$ is NIP.
\end{cor}

Finally, we give a proof of \Cref{NIPnhensexpcon} in the case of residue characteristic $p$:

\begin{cor}\label{KNIPn->KvNIPn}
 Let $(K,v)$ be henselian of residue characteristic $p$. Assume $K$ is \NIPn as a pure field. Then $(K,v)$ is \NIPn as a valued field.
\end{cor}

\begin{proof}Because $K$ is \NIPn and $v$ is henselian, we can apply \Cref{asext}.

If $(K,v)$ is of equicharacteristic $p$, then it is SAMK -- or trivial, but in case the valuation is trivial, there is nothing to prove. Consider $v_K$, the canonical henselian valuation on $K$. We aim to first prove that the residue field $k_{v_K}$ is \NIPn. If $k_{v_K}$ is separably closed, then it is NIP. If $k_{v_K}$ is not separably closed, then $v_K$ is definable; indeed, by \cite[Thm.~3.15]{JK15}, if $k_{v_K}$ is neither separably closed nor real closed and if the absolute Galois group of $K$ is non-universal, then $v_K$ is definable. Here we know that the Galois group is non-universal because a SAMK valued field cannot have Galois extensions of degree $p$-divisible. Since $v_K$ is definable in a \NIPn structure, $k_{v_K}$ is \NIPn. 

Now we prove that $k_v$ is \NIPn. If $v$ is a proper refinement of $v_K$, then $k_v$ is separably closed, hence NIP. If $v$ is a coarsening of $v_K$, then $(k_v,\overline{v_K})$ is SAMK or trivial: indeed, any separable extension of degree $p$-divisible of $k_v$ would lift to a separable extension of same degree for $K$ by henselianity. Hence, since $(k_v,\overline{v_K})$ is SAMK with \NIPn residue, it is \NIPn by \Cref{NIPthm}.

Finally, because $(K,v)$ is SAMK and $k_v$ is \NIPn, we conclude that $(K,v)$ is \NIPn by \Cref{NIPthm}.

Now assume that $(K,v)$ is of mixed characteristic and that $v_p$ is finitely ramified. It is therefore definable in $K$, and so $k_p$ is \NIPn. Now, $(k_p,\overline{v})$ is of equicharacteristic $p$, so it is \NIPn by the argument above. We now apply \Cref{finramtrans} and obtain that $(K,v)$ is \NIPn.

Lastly, if $(K,v)$ is of mixed characteristic and $(k_0,\overline v)$ is AMK, we proceed similarily as is the equicharacteristic $p$ case: if $v$ is a proper refinement of $v_k$, then $k_v$ is separably closed and therefore $(K,v)$ is NIP by \Cref{NIPthm}. If $v$ is a coarsening of $v_K$, we first prove that $k_{v_K}$ is \NIPn, since it is either separably closed, or $v_K$ is definable by \cite[Thm.~3.15]{JK15} -- note that \cite[Obs.~3.16]{JK15} guarantees that the Galois group of $K$ is not universal. Now, we argue as above that $(k_v,\overline{v_K})$ is \NIPn, and then that $(K,v)$ is \NIPn.
\end{proof}

To prove Jahnke's \NIPn Henselian Expansion Conjecture \ref{NIPnhensexpcon}, only the equicharacteristic $0$ case remains. But in this case, there is no reason to believe that the Galois group would be non-universal, and thus no way to define $v_K$, or any other valuation.

 \bibliographystyle{plain}
\bibliography{bib/mtvf.bib,bib/algbooks.bib,bib/mt.bib,bib/vfbooks.bib,bib/mtbooks.bib,bib/misc.bib,bib/alg.bib,bib/itme.bib}
\end{document}